\documentclass[preprint,12pt]{article}

\usepackage[english]{babel}	
\usepackage[utf8x]{inputenc}		
\usepackage{enumerate}
\usepackage{textcomp}                
\usepackage{typearea}
\usepackage{pdfpages}
\usepackage[toc]{appendix}
\usepackage[]{todonotes}
\usepackage{subcaption}
\usepackage{verbatim}
\usepackage{comment}
\usepackage{mathtools}
\mathtoolsset{showonlyrefs}
\usepackage{lscape}

\usepackage{algorithm}
\usepackage{algpseudocode}
\usepackage{algorithmicx}

\usepackage{amsfonts, amsmath, amssymb, amsthm, dsfont}
\usepackage{amsthm}			
\usepackage{thumbpdf}	
\usepackage{natbib}
\usepackage{booktabs}
\usepackage[noblocks]{authblk}
\usepackage{ulem}

\usepackage{pgf, tikz}
\usepackage{tikz-cd}		
\usepackage{bbm}

\theoremstyle{definition}
\newtheorem{definition}{Definition}[]

\newtheorem*{assumption*}{Assumption}

\newtheorem*{condition*}{Condition}

\theoremstyle{plain}
\newtheorem{theorem}[definition]{Theorem}
\newtheorem{proposition}[definition]{Proposition}
\newtheorem{lemma}[definition]{Lemma}

\theoremstyle{remark}
\newtheorem{remark}{Remark}

\usepackage{graphicx}

\newcommand{\N}{\mathbb{N}}
\newcommand{\Z}{\mathbb{Z}}

\newcommand{\EE}{\mathbb{E}}
\newcommand{\R}{\mathbb{R}}

\newcommand{\Le}{\mathcal{L}}

\newcommand{\floor}[1]{\lfloor #1 \rfloor}

\newcommand{\pconv}{\xrightarrow{\mathbb{P}}}

\newcommand{\Cov}{\operatorname{Cov}}

\newcommand{\deq}{\overset{d}{=}}

\newcommand{\wconv}{\Rightarrow}

\newcommand{\tr}{\mathrm{tr}}
\newcommand{\beps}{\boldsymbol{\epsilon}}

\allowdisplaybreaks

\includecomment{proofcomplete}

\title{Strong Gaussian approximations with random multipliers\footnote{Funded by the Federal Ministry of Education and Research (BMBF) and
the Ministry of Culture and Science of the German State of North Rhine-Westphalia
(MKW) under the Excellence Strategy of the Federal Government and the Länder.}}
\author{Fabian Mies\\ Delft University of Technology\footnote{This work was done while the author was employed at RWTH Aachen University, Germany.}}

\begin{document}
	
	\maketitle

	\begin{abstract}
		\noindent
		One reason why standard formulations of the central limit theorems are not applicable in high-dimensional and non-stationary regimes is the lack of a suitable limit object. 
		Instead, suitable distributional approximations can be used, where the approximating object is not constant, but a sequence as well. 
		We extend Gaussian approximation results for the partial sum process by allowing each summand to be multiplied by a data-dependent matrix.
		The results allow for serial dependence of the data, and for high-dimensionality of both the data and the multipliers.
		In the finite-dimensional and locally-stationary setting, we obtain a functional central limit theorem as a direct consequence.
		An application to sequential testing in non-stationary environments is described.
	\end{abstract}

	\section{Introduction}
	
	Donsker's theorem states that for iid, centered, d-variate random vectors $X_t$, the partial sum process $S_n(u) = \frac{1}{\sqrt{n}} \sum_{t=1}^{\lfloor un\rfloor} X_t$, $u\in[0,1]$, converges weakly to $\Sigma^{\frac{1}{2}} B(u)$ in the Skorokhod space $D[0,1]$, where $B(u)$ is a standard d-variate Brownian motion $B(u)$, and $\Sigma=\Cov(X_t)$ is the covariance matrix of the random vectors.
	If the $X_t$ are not iid but stationary and satisfy suitable ergodicity properties, e.g.\ strong mixing, then the same functional central limit theorem (FCLT) holds, but with $\Sigma = \sum_{h=-\infty}^\infty \Cov(X_t, X_{t+h})$ the long-run variance \citep{Billingsley1999}.
	Another deviation from the iid assumption is a potential non-stationarity of the $X_t$. 
	In the latter case, it is difficult to formulate a FCLT because it is in general unclear how to specify a limit random element or distribution. 
	One mathematical trick to solve this challenge is the asymptotic framework of locally stationary time series introduced by \cite{Dahlhaus1997}.
	Here, the non-stationarity is reparametrized from the absolute time $t\in \{1,\ldots, n\}$ to the relative time $\frac{t}{n}\in [0,1]$, such that for $\frac{t}{n}\approx u$, the time series $X_t$ may be well approximated by a stationary time series $X_t(u)$.
	A rather general model formalizing this idea is introduced in \cite{Dahlhaus2017}, see also Section \ref{sec:sequential} below.
	Under suitable conditions, the locally-stationary version of Donsker's Theorem is of the form $S_n(u) \wconv \int_0^u \Sigma^{\frac{1}{2}}(u)\, dB(u)$, where $\Sigma(u)=\sum_{h=-\infty}^\infty \Cov(X_t(u), X_{t+h}(u))$ is the long-run variance of the local approximating time series.
	
	Moving beyond the locally stationary framework, the FCLT breaks down because there is no well-defined limit object. 
	In the same way, the FCLT breaks down in a high-dimensional asymptotic setting, where $d=d_n\to \infty$ as $n\to \infty$. 
	What can be salvaged in these situations is a good approximation of $S_n(u)$ in terms of a Gaussian process $B_n(u)$, say, such that $\sup_{u\in [0,1]} \|S_n(u)-B_n(u)\|\to 0$ in probability as $n\to\infty$. 
	That is, we approximate the sequence $S_n(u)$ by a sequence of Gaussian processes, instead of a single limiting Gaussian process. 
	Approximations of this kind have been known for iid data as strong approximations, or Hungarian couplings \citep{komlos1975,komlos1976}, where optimal rates of approximation have been achieved; see \cite{Zaitsev2013} for an overview. 
	More recent contributions for multivariate and / or dependent data are due to \cite{liu2009, Wu2011,Berkes2014, Karmakar2020,bonnerjee_gaussian_2024}.
	For non-stationary and high-dimensional time series, a corresponding result has been derived in \cite{mies_sequential_2023}; see also \cite{mies_projection_2024}. 
	In this article, we extend the latter approximation in two regards:
	\begin{itemize}
		\item[(i)] We exploit a stronger covariance bound to improve the approximation rates of \cite{mies_sequential_2023}.
		\item[(ii)] We extend the results to partial sums with data-dependent multipliers, i.e.\ $\widetilde{S}_n(u) = \frac{1}{\sqrt{n}}\sum_{t=1}^{\lfloor un \rfloor} \widehat{g}_{t}X_t$ for random matrices $\widehat{g}_t\in \R^{m\times d}$, with potentially growing dimension $d\to \infty$, $m\to\infty$.
	\end{itemize}	
	
	\paragraph{Notation}
	For a vector $x\in\R^d$, we denote by $\|x\|=(\sum_{i=1}^d x_i^2)^{1/2}$ the Euclidean norm, and for a matrix $A\in \R^{d\times d}$, $A^T$ denotes its transpose, $\tr(A)$ its trace, $\|A\|_{\tr} = \tr((A^T A)^{\frac{1}{2}})$ its trace norm.
	For any matrix $A\in\R^{m\times d}$, we denote by $\|A\|=\|A\|_F = \sqrt{\tr(A^TA)}$ the Frobenius norm. 
	We use $\wconv$ to denote weak convergence of probability measures resp.\ random elements.
	The $\Le_q$ norm of a random vector $X$ is denoted as $\|X\|_{\Le_q} = (\EE \|X\|^q)^{1/q}$, i.e.\ we always use the Euclidean base-norm. 
	For two symmetric matrices $A,B$, the notation $A\succeq B$ means that $A-B$ is positive semidefinite.
	
	\section{Sequential approximations for time series}\label{sec:sequential}
	
	A general class of non-stationary and potentially nonlinear time series may be formulated as follows.	
	For iid random seeds $\epsilon_i \sim U(0,1)$, and non-random functions $G_{t,n}:\R^\infty \to \R^d$, $t=1,\ldots, n$, define the sequence of non-stationary time series $\{X_{t,n}, t=1,\ldots,n\}$ with $X_{t,n}\in\R^d$ as 
	\begin{align*}
		X_{t,n} &= G_{t,n}(\beps_t), \quad t=1,\ldots, n, \\
		\beps_t &= (\epsilon_t, \epsilon_{t-1},\ldots) \in\R^\infty.
	\end{align*}
	Throughout this section, we assume that $\EE(X_{t,n})=0$.
	For an independent copy $\tilde{\epsilon}_i\sim U(0,1)$, define furthermore
	\begin{align*}
		\tilde{\beps}_{t,h} &= (\epsilon_t, \ldots, \epsilon_{t-h+1}, \tilde{\epsilon}_{t-h},\epsilon_{t-h-1}\ldots) \in \R^\infty, \\
		\beps_{t,h} &= (\epsilon_t, \ldots, \epsilon_{t-h+1}, \tilde{\epsilon}_{t-h},\tilde{\epsilon}_{t-h-1}\ldots) \in \R^\infty.
	\end{align*}
	To express the mixing-type behavior of $X_{t,n}$, we employ the physical dependence measure introduced by \cite{Wu2005}, and defined as
	\begin{align*}
		\delta_{n}(h) := \max_{t=1,\ldots, n} \|G_{t,n}(\beps_0) - G_{t,n}(\tilde{\beps}_{0,h}) \|_{\Le_q}.
	\end{align*}
	for some $q\geq 2$. 
	The dependence of $\delta_n$ on $q$ will be implicit in the remainder of this article.
	
	\begin{remark}
		For a nonstationary linear process of the form $X_{t,n} = \sum_{i=0}^\infty a_{t,n,i}\eta_{t-i}$ with iid random vectors $\eta_{t}\in\R^d$ and matrices $a_{t,n,i}\in\R^{d\times d}$, we may write $\eta_t=h(\epsilon_t)$ for a measurable mapping $h:\R\to\R^d$. 
		The physical dependence measure may then be computed as $\delta_n(h)=\|\eta_1\|_{\Le_q}\sup_{t} \|\alpha_{t,n,h}\|_{\text{op}}$, where $\|A\|_{\text{op}}$ denotes the Euclidean operator norm of a matrix $A\in\R^{d\times d}$.
		Upper bounds on the physical dependence measure may also be established for a broader class of time-varying ARMA processes as in \cite{subbarao1970, grenier1983, moulines2005, Dahlhaus2009a, Giraud2015}.
		Here, $\delta_n(h)\to 0$ polynomially as $h\to \infty$ if the roots of the auto-regressive poolynomial are bounded away from one \citep[Example 1]{mies_functional_2023}.
	\end{remark}
	
	We first impose the following set of assumptions, for some $\Gamma_n\geq 1$ and $\beta>1$:
	\begin{align}
		\|G_{1,n}(\beps_0)\|_{\Le_2} + \sum_{t=2}^n \|G_{t,n}(\beps_0) - G_{t-1,n}(\beps_0)\|_{\Le_2}  &\leq \Theta_n \Gamma_n, \tag{A.1} \label{ass:A1} \\
		\delta_n(h) &\leq \Theta_n (h+1)^{-\beta}, \tag{A.2} \label{ass:A2}.
	\end{align}
	As a first result, we find that the Gaussian approximation result presented in \cite{mies_sequential_2023} can be slightly improved in terms of the dependence on the exponent $\beta$. 
	In particular, the following property has not been fully exploited in the proofs therein.
	
	\begin{proposition}\label{prop:autocov}
		Assumption \eqref{ass:A2} implies that $\|\Cov(X_{t,n}, X_{s,n})\|_{\tr} \leq \Theta_n^2 C_\beta (|s-t|+1)^{-\beta}$.
	\end{proposition}

	To state the Gaussian approximation result, define
	\begin{align*}
		\xi(q,\beta) = 
		\begin{cases}
			\frac{q-2}{6q-4}, 
			& \beta\geq \frac{3}{2},\, \beta> \frac{2q}{q+2}, \\
			\frac{\beta-1}{4\beta-2}, 
			& \beta\geq \frac{3}{2},\, \beta \leq \frac{2q}{q+2}, \\
			\frac{(\beta-1)(q-2)}{4q\beta-3q-2}, 
			& \beta<\frac{3}{2},\, \beta > \frac{2q}{q+2}, \\
			\frac{(\beta-1)^2}{2\beta^2-1-\beta}, 
			& \beta<\frac{3}{2},\, \beta\leq \frac{2q}{q+2},
		\end{cases}
		\qquad q>2,\, \beta>1.
	\end{align*}
	
	\begin{theorem}\label{thm:Gauss-ts}
		Let $X_{t,n}=G_{t,n}(\beps_t)$ with $\EE (X_{t,n}) = 0$ satisfy \eqref{ass:A1} and \eqref{ass:A2}, for some $q>2$, and suppose $d\leq c n$ for some $c>0$. 
		Suppose furthermore that $\beta>1$, which implies that the local long-run covariance matrix $\Sigma_{t,n} = \sum_{h=-\infty}^\infty \Cov(G_{t,n}(\beps_0),G_{t,n}(\beps_h))$ is well-defined.
		Then, on a potentially different probability space, there exist random vectors $(X_{t,n}')_{t=1}^n\deq (X_{t,n})_{t=1}^n$ and independent, mean zero, Gaussian random vectors $Y_t^*\sim\mathcal{N}(0,\Sigma_{t,n})$ such that
		\begin{align}
			\left(\EE \max_{k\leq n} \left\|\frac{1}{\sqrt{n}}\sum_{t=1}^k  (X_{t,n}' - Y_{t,n}^*)   \right\|^2 \right)^\frac{1}{2} 
			&\leq C\Theta_n  \Gamma_n^{ \frac{\beta-1}{2\beta}} \sqrt{\log(n)} \left( \frac{d}{n} \right)^{\xi(q,\beta)}. 
			\label{eqn:Gauss-ts-2}
		\end{align}
	\end{theorem}
	
	\begin{remark}
		The essential improvement of Theorem \ref{thm:Gauss-ts} over Theorem 3.1(ii) of \cite{mies_sequential_2023} is that the new formulation only requires $\beta>1$, instead of $\beta>2$. 
		This is possible by exploiting Proposition \ref{prop:autocov} throughout the proof. 
	\end{remark}

	\section{Strong approximation with random multipliers}
	
	In this section, we extend the Gaussian approximation result to also allow for data-dependent multipliers. 
	That is, we consider the partial sums $\frac{1}{\sqrt{n}} \sum_{t=1}^n \hat{g}_{t,n} X_{t,n}$ for random matrices $\hat{g}_{t,n}\in \R^{m\times d}$, with potentially growing dimensions $d=d_n\to \infty$, $m=m_n\to \infty$.
	The matrices $\hat{g}_{t,n}$ are assumed to converge to some non-random matrix $g_{t,n}$ in the sense that
	\begin{align*}
		\Lambda_n :=\sqrt{\sum_{t=1}^n \|\hat{g}_{t,n} - g_{t,n}\|_{\Le_2}^2} \, = o(\sqrt{n}).
	\end{align*}
	Moreover, the latter sequences of matrices should admit some regularity in time, which we express in terms of the quantities 
	\begin{align*}
		\Psi_n :=\sum_{t=2}^n \|g_{t,n} - g_{t-1,n}\|,\qquad \Phi_n = \max_{t=1,\ldots,n} \|g_{t,n}\|.
	\end{align*}
	In all the terms $\Lambda_n$, $\Psi_n$, and $\Phi_n$, the dependence on the dimension $d_n$ and $m_n$ is implicit in the norms.
	Besides these technical restrictions, we impose the important qualitative restriction that, for some sequence $L_n\in \N$ which is not necessarily diverging, 
	\begin{align}
		\text{$\widehat{g}_{t,n}$ is $\beps_{t- L_n}$-measurable.}\tag{A.M} \label{ass:AX}
	\end{align}
	Thus, if $\widehat{g}_{t,n}$ is based on an estimator, then it may only use the data up to time $t-L_n$. 
	A similar measurability assumption was imposed in \cite{mies_functional_2023}.
	To express the relation between the delay $L_n$ and the serial dependence of the time series, define the term
	\begin{align*}
		\Xi_n := \sum_{h=L_n}^\infty \delta_n(h).
	\end{align*}
	\begin{remark}
		Note that \eqref{ass:A2} implies that $\Xi_n \leq \Phi_n L_n^{1-\beta}$, but in some situations, it might be possible to derive sharper bounds on $\Xi_n$.
		In particular, if the $X_{t,n}$ are indeed independent, then $\Xi_n=0$ for any sequence $L_n\geq 1$, thus greatly reducing the measurability conditions on $\widehat{g}_{t,n}$. 
		In the latter case, \eqref{ass:A2} also holds for all $\beta>0$.
	\end{remark}
	The local long run covariance matrix of the process $g_{t,n} X_{t,n}$ is given by
	\begin{align*}
		\Sigma_{t,n}
		&= \sum_{h=-\infty}^\infty g_{t,n} \Cov[G_{t,n}(\beps_0), G_{t,n}(\beps_h)] g_{t,n}^T.
	\end{align*}
	Note that $\Sigma_{t,n}$ is finite for $\beta>1$ by virtue of Proposition \ref{prop:autocov}.
	
	\begin{theorem}\label{thm:multiplier}
		Let Assumptions \eqref{ass:A1} and \eqref{ass:A2} hold, and suppose that $\hat{g}_{t,n}$ is $\beps_{t-L_n}$-measurable for some $L_n\geq 1$.
		Then, upon potentially enlarging the probability space, there exist independent Gaussian random vectors $Y_{t,n}\sim\mathcal{N}(0,\Sigma_{t,n})$, such that
		\begin{align*}
			&\quad \left| \max_{k=1,\ldots,n} \frac{1}{\sqrt{n}} \sum_{t=1}^k [\hat{g}_{t,n} X_{t,n} - Y_{t,n}] \right| \\
			&= \mathcal{O}_P\left( n^{-\frac{1}{2}} \Lambda_n \Theta_n + \Lambda_n \Xi_n + n^{\frac{1}{q}-\frac{1}{2}} L_n \Phi_n\right) \\
			&\quad + \mathcal{O}_P\left( \Phi_n \Theta_n  \left( \Phi_n\Theta_n \Gamma_n + \Psi_n\Theta_n \right)^{\frac{\beta-1}{2\beta}} \sqrt{\log(n)} \left( \tfrac{m}{n} \right)^{\xi(q,\beta)} \right).
		\end{align*}
		This result is also applicable if $m=m_n\to\infty$.
	\end{theorem}
	
	Multiplier central limit theorems are often used in the literature to study the behavior of bootstrap methods, see e.g.\ \cite[Sec.~2.9]{vandervaart1996}.
	There, the process $X_{t,n}$ and the multipliers are independent, which is not the case in our setting. 
	Instead, Assumption \eqref{ass:A4} implies that the $\hat{g}_{t,n}$ are almost deterministic, up to some estimation error. 
	The contribution of Theorem \ref{thm:multiplier} is to take this estimation error explicitly into account.
	Compared to Theorem \ref{thm:Gauss-ts}, it allows for greater flexibility in the decay of the physical dependence measure $\delta_n(h)$ as expressed via the quantity $\Xi_n$, rather than assuming a polynomial decay.

	\section{Functional central limit theorem}
	
	In the classical finite-dimensional setting where $d$ and $m$ are constant, the previous results also allow for the formulation of a functional central limit theorem. 
	We embed the idea of local stationarity in our framework by considering a local asymptotic mapping $G_u:\R^\infty \to \R^d$, $u\in [0,1]$, such that, as $n\to\infty$,
	\begin{align}
		\int_0^1\|G_{\lfloor vn\rfloor,n}(\beps_0) - G_{v}(\beps_0)\|_{\Le_2}\, dv &\to 0.\tag{A.3}\label{ass:A3} 
	\end{align}
	Then $X_t(u)=G_u(\beps_t)$ is the local stationary approximation in the sense of \cite{Dahlhaus1997}.
	Analogously, we suppose that there exists another class of matrices $g_u\in \R^{m\times d}$, $u\in[0,1]$, such that $\sup_u \|g_u\| < \infty$ and, as $n\to\infty$,
	\begin{align}
		\int_0^1 \|g_{\lfloor vn\rfloor,n} - g_v\|\, dv &\to 0. \tag{A.4} \label{ass:A4}
	\end{align}
	Since we are now working with fixed dimensions $d$ and $m$, the choice of the matrix norm is irrelevant.
	To formulate the functional central limit theorem, introduce the asymptotic local long-run variance matrix
	\begin{align*}
		\Sigma_u &= \sum_{h=-\infty}^\infty g_{u} \Cov[G_{u}(\beps_0), G_{u}(\beps_h)] g_{u}^T.
	\end{align*}
	
	\begin{lemma}\label{lem:FCLT}
		Suppose that \eqref{ass:A2}, \eqref{ass:A3}, \eqref{ass:A4} hold, and $\Theta_n,\Phi_n=\mathcal{O}(1)$. 
		For independent Gaussian random vectors $Y_{t,n} \sim \mathcal{N}(0, \Sigma_{t,n})$, we have the following weak convergence in the Skorokhod space $D[0,1]$:
		\[
		\frac{1}{\sqrt{n}} \sum_{t=1}^{\lfloor un\rfloor} Y_{t,n} \wconv \int_0^u \Sigma_v^{\frac{1}{2}} \, dW_v.
		\]
	\end{lemma}

	Combining Theorem \ref{thm:multiplier} and Lemma \ref{lem:FCLT}, we immediately obtain the main result of this paper.
	
	\begin{theorem}\label{thm:multiplier-FCLT}
		Suppose that $\hat{g}_{t,n}$ is $\beps_{t-L}$-measurable for $L=L_n$, and let \eqref{ass:A1}, \eqref{ass:A2}, \eqref{ass:A3}, \eqref{ass:A4} hold with $\Theta_n, \Phi_n=\mathcal{O}(1)$ such that
		\begin{align*}
			\Lambda_n\left( n^{-\frac{1}{2}}  +  \Xi_n\right) \to 0, \qquad L_n = o\left(n^{\frac{1}{2}-\frac{1}{q}}\right), \\
			\left( \Gamma_n + \Psi_n \right)^{\frac{\beta-1}{2\beta}} \sqrt{\log(n)} n^{-\xi(q,\beta)} \to 0.
		\end{align*}
		Then
		\begin{align*}
			\frac{1}{\sqrt{n}} \sum_{t=1}^{\floor{un}} \hat{g}_{t,n} [X_{t,n}-\EE(X_{t,n})] \;\wconv\; \int_0^u \Sigma_v^\frac{1}{2}\, dW_v,
		\end{align*}
		where weak convergence holds in the Skorokhod space $D[0,1]$.
	\end{theorem}
	
	\section{Example: studentized partial sums}
	
	In the sequel, we outline a statistical application of the FCLT with data-dependent multipliers.
	Consider a heteroskedastic sequence of independent $d$-variate random vectors $X_{t,n} = \Sigma(t/n)^{\frac{1}{2}} \eta_t$, for iid random vectors $\eta_t$ with $\Cov(\eta_t) = I_{d\times d}$, $\EE(\eta_t)=0$, $\EE \|\eta_t\|^q<\infty$ for some $q>2$, and a Lipschitz-continuous function $\Sigma:[0,1] \to \R^{d\times d}$. Suppose that each $\Sigma(u)$ is symmetric positive definite, and that $\Sigma(u)$ is uniformly elliptic: there exists some small $c>0$ such that $\Sigma(u)\succeq cI_{d\times d}$, i.e.\ $\Sigma(u)-cI_{d\times d}$ is positive semidefinite for all $u\in[0,1]$.
	The FCLT for locally-stationary data, e.g.\ Theorem \ref{thm:multiplier-FCLT} with $\widehat{g}_{t,n}=I_{d\times d}$, yields that
	\begin{align}
		S_n(u)=\frac{1}{\sqrt{n}} \sum_{t=1}^{\lfloor un\rfloor} X_{t,n} \wconv \int_0^u \Sigma(v)^{\frac{1}{2}} dW_v. \label{eq:not-studentized}
	\end{align}
	Although technically a special case of the previous theorem, the latter convergence is rather standard in the literature and does not exploit our new results. 
	Instead, to showcase the application of Theorem \ref{thm:multiplier-FCLT}, we estimate $\Sigma(u)$ locally by the window average
	\begin{align*}
		\widehat{\Sigma}(t/n) = \frac{1}{k_n} \sum_{j=1}^{k_n} X_{t-j} X_{t-j}^T, \qquad t=k_n+1,\ldots, n.
	\end{align*}
	For $t\leq k_n$, let $\widehat{\Sigma}(t/n)=0$.
	To enforce the uniform ellipticity, we define
	\begin{align*}
		\widetilde{\Sigma}(t/n) =\begin{cases}
			\widehat{\Sigma}(t/n), & \text{if } \widehat{\Sigma}(t/n) \succeq cI_{d\times d}, \\
			c I_{d\times d}, & \text{otherwise}.
		\end{cases} 
	\end{align*}
	We consider the studentized partial sum process
	\begin{align*}
		S_n^*(u) = \frac{1}{\sqrt{n}} \sum_{t=k_n+1}^{\lfloor un \rfloor} \widetilde{\Sigma}(t/n)^{-\frac{1}{2}} X_{t,n}. 
	\end{align*}
	
	\begin{theorem}\label{thm:studentized}
		For $k_n = \lfloor n^{\frac{2}{3}}\rfloor$, the studentized partial sum process converges weakly in the Skorokhod space to a standard Brownian motion, 
		\begin{align}
			S_n^*(u) \wconv \int_0^u \, dW_v. \label{eq:studentized}
		\end{align}
	\end{theorem}

	Both limits \eqref{eq:not-studentized} and \eqref{eq:studentized} can be used to  test the null hypothesis $H_0: \EE(X_{t,n})=0$ via the test statistic 
	\begin{align*}
		T_n^* &= \sup_{u\in [0,1]} \|S_n^*(u)\| \wconv \sup_{u\in [0,1]} \left\|\int_0^u dW_v\right\|, \\
		T_n &= \sup_{u\in [0,1]} \|S_n(u)\| \wconv \sup_{u\in [0,1]} \left\|\int_0^u \Sigma(v)^{\frac{1}{2}} dW_v\right\|.
	\end{align*}
	The supremum may be evaluated exactly as $S_n$ and $S_n^*$ are step functions. 
	Because $T_n^*$ weights the observations according to their inverse variance, we expect this statistic to be more powerful. This is supported by simulation results reported in Table \ref{tab:power}, for sample size $n=10000$. 
	Here, the data is generated as $X_{t,n} = \sigma(\frac{t}{n}) (Z_t-1) + \mu$ with $Z_t \sim \text{Exp}(1)$, and $\sigma(u) = 1.2 + \sin(6\pi u)$.
	It is found that both the size distortion of $T_n^*$ is slightly bigger, but its power is much higher than for the test statistic $T_n$. However, for smaller sample sizes ($n=1000$), we observed that $T_n^*$ suffers from a more pronounced upward size distortion (test size $0.097$).
	
	A second important distinction is that $T_n^*$ enables a sequential testing procedure. 
	Since the critical value $c_\alpha^*$ of $\sup_u \|W_u\|$ does not depend on the unknown function $\Sigma$, we can reject $H_0$ as soon as $\|S_n^*(u)\|>c_\alpha^*$. 
	Although the non-studentized test statistic $\|S_n(u)\|$ may also be computed sequentially, the corresponding critical value $c_\alpha(\Sigma)$ will depend on the whole function $\Sigma(u)$. However, at time $t$, the best we can achieve is to estimate $\Sigma(u)$ for $u\in [0, \frac{t}{n}]$.
	
	\begin{table}
		\centering
		\begin{tabular}{r|ccccccccccc}
			\toprule
			$\mu$ & 0 & .002 & .004 & .006& .008 & .010 & .012 & .014 & .016 & .018 & .020 \\ \midrule
			Power of $T_n$ & .048 & .053 & .054 & .060 & .078 & .094 & .119 & .151 & .184 &  .213 & .254 \\
			Power of $T_n^*$ & .066 & .079 & .097 & .134 & .194 & .261 & .362 & .456 & .559 & .656 & .743 \\ \bottomrule
		\end{tabular}
		\caption{Power of the test statistics $T_n$ and $T_n^*$ for sample size $n=10000$. Reported results are based on $10^4$ simulations each.}
		\label{tab:power}
	\end{table}
	
	The idea of local studentization can also be extended to non-stationary time series. In this situation, we need to replace $\widehat{\Sigma}(t/n)$ by a suitable estimator of the long-run variance matrix. 
	The technical details are beyond the scope of this article.

	\section*{Proofs}
	
	\subsection*{Proof of Proposition \ref{prop:autocov}}
	Observe that the random vectors $Y_{t,n,i} = \EE(X_{t,n}|\beps_i) - \EE(X_{t,n}|\beps_{i-1})$, $i\in\Z$, are martingale differences, with
	\begin{align*}
		\|Y_{t,n,i}\|_{\Le_q} &\leq \begin{cases} \|G_{t,n}(\beps_t) - G_{t,n}(\beps_{t,t-i})\|_{\Le_q} \leq \Theta_n (t-i+1)^{-\beta},& \text{ if } i\leq t, \\
			0,&  \text{ if } i>t. 
		\end{cases}
	\end{align*}
	The analogous bound holds for $Y_{s,n,i} = \EE(X_{s,n}|\beps_i) - \EE(X_{s,n}|\beps_{i-1})$. 
	Thus,
	\begin{align*}
		\|\Cov(X_{t,n}, X_{s,n})\|_\tr &= \left\| \sum_{i\in \Z} \Cov(Y_{t,n,i}, Y_{s,n,i}) \right\|_\tr \\
		&\leq \sum_{i\leq (s\wedge t)} \EE \|Y_{t,n,i}Y_{s,n,i}^T\|_\tr \\
		&\overset{(\dagger)}{=} \sum_{i\leq (s\wedge t)} \EE \left( \|Y_{t,n,i}\| \cdot \|Y_{s,n,i}\| \right) \\
		&\leq \sum_{i\leq (s\wedge t)} \|Y_{t,n,i}\|_{\Le_2} \cdot \|Y_{s,n,i}\|_{\Le_2} \\
		&\leq \Theta_n^2 \sum_{k=0}^\infty (k+1)^{-\beta} (k+|s-t|+1)^{-\beta} \\
		&\leq \Theta_n^2 C_\beta (|s-t|+1)^{-\beta}.
	\end{align*}
	At the step $(\dagger)$, we use that $\|x y^T\|_\tr = \|x\| \cdot \|y\|$, see for example \cite[Lemma 6.3]{mies_sequential_2023}.
	
	\subsection*{Proof of Theorem \ref{thm:Gauss-ts}}
	We show how to improve various steps in the proof of Theorem 3.1 of \cite{mies_sequential_2023}. 
	For consistency with the latter reference, in this proof, $L$ denotes a block length parameter distinct from the lag sequence $L_n$ of the main text. 
	
	By using the improved inequality $\|\Cov(G_{t,n}(\beps_r), G_{t,n}(\beps_s))\|_F \leq C \Theta^2 (|r-s|+1)^{-\beta}$ from Proposition \ref{prop:autocov}, we may improve (22) therein to
	\begin{align*}
		\left\| \Sigma_{t,n} - \Sigma_{t,n}^l \right\|_\tr  &\leq C \Theta_n^2 L^{(1-\beta)\vee (-1)},
	\end{align*}
	and (23) therein to
	\begin{align*}
		\|\Sigma_{t,n}^l - \Cov(Y_{t,n}')\|_{\tr} 
		&\leq C\Theta_n^2 \left[ \frac{1}{\Theta_n}\sum_{s=t_{l-1}+1}^{t_l} \|G_{s,n}(\beps_0)-G_{s-1,n}(\beps_0)\|_{\Le_2} \right]^{1-\frac{1}{\beta}}.
	\end{align*}
	Hence, we can improve (26) therein to
	\begin{align*}
		&\quad \left( \EE \max_{k\leq n} \left\| \sum_{t=1}^k Y_{t,n}^* - Y_{t,n}'\right\|^2 \right)^\frac{1}{2} \\
		&\leq C \Theta_n n^\frac{1}{2} L^{ \frac{1-\beta}{2} \vee (-\frac{1}{2})  } + C\Theta_n L^\frac{1}{2} \Gamma_n^{\frac{\beta-1}{2\beta}} M^{\frac{1}{2\beta}}.
	\end{align*}
	That is, all occurrences of $\beta$ in (26) therein are replaced by $\beta+1$. 
	As in (27) therein, we obtain for any $L\leq c\frac{n}{d}$,
	\begin{align}
		&\left( \EE \max_{k\leq n} \left\|\frac{1}{\sqrt{n}}\sum_{t=1}^k (X_{t,n}' -Y_{t,n}^*) \right\|^2 \right)^\frac{1}{2} \nonumber \\
		&\leq C \Theta_n \sqrt{\log(n)} \left( \frac{d L}{n} \right)^{\frac{1}{4}-\frac{1}{2q}} + C\Theta \left( L^{1-\beta} + L^{-\frac{1}{2}} \right) \nonumber \\
		&\qquad+C \Theta_n  L^{ \frac{1-\beta}{2} \vee (-\frac{1}{2})  } + C\Theta_n n^{-\frac{1}{2}} L^\frac{1}{2} \Gamma_n^{\frac{\beta-1}{2\beta}} M^{\frac{1}{2\beta}} \nonumber \\
		&\leq C\Theta_n \Gamma_n^{ \frac{\beta-1}{2\beta}} \sqrt{\log(n)} \left\{  \left( \frac{d L}{n} \right)^{\frac{q-2}{4q}} +  L^{1-\beta} + L^{-\frac{1}{2}}  + \left(\frac{L}{n} \right)^{\frac{\beta-1}{2\beta}}\right\}  \nonumber \\
		&\leq \begin{cases}
			C\Theta_n \Gamma_n^{ \frac{\beta-1}{2\beta}} \sqrt{\log(n)} \left\{  \left( \frac{d L}{n} \right)^{\frac{q-2}{4q}} + L^{-\frac{1}{2}}\right\}, 
			& \beta\geq \frac{3}{2},\, \beta> \frac{2q}{q+2}, \\
			C\Theta_n \Gamma_n^{ \frac{\beta-1}{2\beta}} \sqrt{\log(n)} \left\{  \left( \frac{d L}{n} \right)^{\frac{\beta-1}{2\beta}} + L^{-\frac{1}{2}}\right\}, 
			& \beta\geq \frac{3}{2},\, \beta\leq \frac{2q}{q+2}, \\
			C\Theta_n \Gamma_n^{ \frac{\beta-1}{2\beta}} \sqrt{\log(n)} \left\{  \left( \frac{d L}{n} \right)^{\frac{q-2}{4q}} + L^{1-\beta} \right\}, 
			& \beta<\frac{3}{2},\, \beta> \frac{2q}{q+2}, \\
			C\Theta_n \Gamma_n^{ \frac{\beta-1}{2\beta}} \sqrt{\log(n)} \left\{  \left( \frac{d L}{n} \right)^{\frac{\beta-1}{2\beta}} + L^{1-\beta} \right\}, 
			& \beta<\frac{3}{2},\, \beta\leq \frac{2q}{q+2}.
		\end{cases}
		\label{eqn:approx-rate-r2-full}
	\end{align}
	The value of $L$ has been left unspecified and may now be chosen. 
	In the three cases, the rate-optimal choice $L^*$ of $L$ is, respectively,
	\begin{align*}
		L^* \asymp  \begin{cases}
			(\tfrac{n}{d})^{\frac{q-2}{3q-2}}, 
			& \beta\geq \frac{3}{2},\, \beta> \frac{2q}{q+2}, \\
			(\tfrac{n}{d})^{\frac{\beta-1}{2\beta-1}}, 
			& \beta\geq \frac{3}{2},\, \beta \leq \frac{2q}{q+2}, \\
			(\tfrac{n}{d})^{ \frac{q-2}{4q\beta-3q-2} }, 
			& \beta<\frac{3}{2},\, \beta > \frac{2q}{q+2}, \\
			(\tfrac{n}{d})^{ \frac{\beta-1}{2\beta^2-1-\beta} }, 
			& \beta<\frac{3}{2},\, \beta\leq \frac{2q}{q+2}.
		\end{cases}
	\end{align*}
	These choices achieve the rates \eqref{eqn:Gauss-ts-2}.

	\subsection*{Proof of Theorem \ref{thm:multiplier}}
	Denote $L=L_n$.
	
	\underline{De-randomizing the multiplier:}
	Assume w.l.o.g.\ that $X_{t,n}$ is centered, and denote $e_t = \hat{g}_{t,n}-g_{t,n}$.
	Decompose
	\begin{align}
		\sum_{t=1}^k e_t X_{t,n} 
		&= \sum_{j=0}^{\lceil \frac{n}{L}\rceil-1} \sum_{t=jL +1}^{(j+1)L \wedge k} e_t X_{t,n} \nonumber \\
		&= \sum_{j=0}^{\lceil \frac{n}{L}\rceil-1} \sum_{t=jL +1}^{(j+1)L \wedge k} e_t \widetilde{X}_{t,n} 
		+  \sum_{j=0}^{\lceil \frac{n}{L}\rceil-1} \sum_{t=jL +1}^{(j+1)L \wedge k} e_t [X_{t,n}-\widetilde{X}_{t,n}], \label{eqn:multiplier-1} \\
		\text{for} \quad \widetilde{X}_{t,n} &= G_{t,n}(\tilde{\beps}^j_{t,t-jL}), \qquad t=jL+1,\ldots, (j+1)L. \nonumber
	\end{align}
	Here, $\beps_{t,h} = (\epsilon_t, \ldots, \epsilon_{t-h+1}, \tilde{\epsilon}_{t-h}^j,\tilde{\epsilon}_{t-h-1}^j,\ldots)$ where $\tilde{\epsilon}_t^j\sim U(0,1)$ are independent for all $t$ and all $j$.
	The decomposition is chosen such that the terms $\sum_{t=jL+1}^{(j+1)L} e_t \widetilde{X}_{t,n}$ are martingale differences, and $e_t$ and $\widetilde{X}_{t,n}$ are independent.
	Hence
	\begin{align*}
		\EE\max_{r=1,\ldots, \lfloor \frac{n}{L}\rfloor -1} \left\| \sum_{j=0}^{r} \sum_{t=jL +1}^{(j+1)L } e_t \widetilde{X}_{t,n} \right\|^2 
		&= \sum_{j=0}^{\lfloor\frac{n}{L}\rfloor -1}  \EE\left\| \sum_{t=jL +1}^{(j+1)L } e_t \widetilde{X}_{t,n} \right\|^2  \\
		&= \sum_{j=0}^{\lfloor\frac{n}{L}\rfloor -1} \sum_{s,t=jL +1}^{(j+1)L } \tr \left[ \EE \left( e_t \widetilde{X}_{t,n} \widetilde{X}_{s,n}^T e_s^T \right) \right] \\
		&= \sum_{j=0}^{\lfloor\frac{n}{L}\rfloor -1} \sum_{s,t=jL +1}^{(j+1)L }  \EE \left(\tr \left[  e_t \widetilde{X}_{t,n} \widetilde{X}_{s,n}^T e_s^T \right] \right) \\
		&= \sum_{j=0}^{\lfloor\frac{n}{L}\rfloor -1} \sum_{s,t=jL +1}^{(j+1)L }  \EE \left(\tr \left[ e_s^T e_t \widetilde{X}_{t,n} \widetilde{X}_{s,n}^T \right] \right) \\
		&= \sum_{j=0}^{\lfloor\frac{n}{L}\rfloor -1} \sum_{s,t=jL +1}^{(j+1)L } \tr\left[\EE(e_s^T e_t) \EE\left(\widetilde{X}_{t,n}\widetilde{X}_{s,n}^T\right)\right] \\
		&\leq \sum_{j=0}^{\lfloor\frac{n}{L}\rfloor -1} \sum_{s,t=jL +1}^{(j+1)L } \EE \|e_s^T e_t\|_F \cdot  \left\|\EE\left(\widetilde{X}_{t,n}\widetilde{X}_{s,n}^T\right) \right\|_F \\
		&\leq \Theta_n^2 \sum_{j=0}^{\lfloor\frac{n}{L}\rfloor -1} \sum_{s,t=jL +1}^{(j+1)L } \|e_s\|_{\Le_2} \|e_t\|_{\Le_2} (|s-t|+1)^{-\beta} \\
		&\leq  \Theta_n^2 \sum_{j=0}^{\lfloor\frac{n}{L}\rfloor -1} \sum_{s,t=jL +1}^{(j+1)L } (\|e_s\|^2_{\Le_2} +\|e_t\|^2_{\Le_2}) (|s-t|+1)^{-\beta}\\
		&=  \Theta_n^2 \sum_{j=0}^{\lfloor\frac{n}{L}\rfloor -1} \sum_{s,t=jL +1}^{(j+1)L } 2\|e_s\|^2_{\Le_2} (|s-t|+1)^{-\beta} \\
		&\leq   C_\beta \Theta_n^2 \sum_{j=0}^{\lfloor\frac{n}{L}\rfloor -1} \sum_{s,t=jL +1}^{(j+1)L } \|e_s\|^2_{\Le_2} \qquad \text{(because $\beta>1$)}\\
		&= C_\beta \Theta_n^2 \Lambda_n^2.
	\end{align*}
	Moreover, the block discretization error may be bounded uniformly in $k$ as
	\begin{align*}
		\max_{k=1,\ldots, n} \max_{j=0,\ldots,\lceil \frac{n}{L} \rceil -1} \sum_{t=(jL+1)\vee k}^{(j+1)L\wedge k} e_t \widetilde{X}_{t,n}
		&\leq \Phi_n L_n \max_{t=1,\ldots, n} |\widetilde{X}_{t,n}| 
		= \mathcal{O}_P\left(\Phi_n L_n n^{\frac{1}{q}}\right).
	\end{align*}
	
	Regarding the second term in \eqref{eqn:multiplier-1}, we find that
	\begin{align*}
		\max_{k=1,\ldots,n} \left|\sum_{j=0}^{\lceil \frac{n}{B}\rceil-1} \sum_{t=jB +1}^{(j+1)B \wedge k} e_t [X_{t,n}-\widetilde{X}_{t,n}] \right|
		&\leq  \left[ \sum_{t=1}^n e_t^2 \right]^{\frac{1}{2}} \left[ \sum_{j=0}^{\lceil \frac{n}{B}\rceil-1} \sum_{t=jB +1}^{(j+1)B \wedge n}  [X_{t,n}-\widetilde{X}_{t,n}]^2 \right]^{\frac{1}{2}} \\
		&= \mathcal{O}_P \left(  \Lambda_n  \right) \mathcal{O}_P\left( \sqrt{n} \Xi_n\right) 
	\end{align*}
	We have thus established that
	\begin{align*}
		\frac{1}{\sqrt{n}}\sum_{t=1}^k \hat{g}_{t,n} X_{t,n}  = \frac{1}{\sqrt{n}}\sum_{t=1}^k g_{t,n} X_{t,n} + \mathcal{O}_P\left( n^{-\frac{1}{2}} \Lambda_n \Theta_n + \Lambda_n \Xi_n + n^{\frac{1}{q}-\frac{1}{2}} L_n \Phi_n\right),
	\end{align*}
	uniformly in $k=1,\ldots,n$.
	
	\underline{Sequential Gaussian approximation:}
	Denote $Z_{t,n}=H_{t,n} = g_{t,n} G_{t,n}$, such that $g_{t,n} X_{t,n} = H_{t,n}(\beps_t)$, and long-run covariance matrix $\Sigma_{t,n} = \sum_{h=\infty}^\infty \Cov(H_{t,n}(\beps_0), H_{t,n}(\beps_h))$.
	Because $g_{t,n}$ is bounded, this time series satisfies $\|H_{t,n}(\beps_0) - H_{t,n}(\tilde{\beps}_{0,h})\|_{\Le_q} \leq \Phi_n \Theta_n (h+1)^{-\beta}$.
	Moreover,
	\begin{align*}
		\sum_{t=2}^n \left\| H_{t,n}(\beps_0) - H_{t-1,n}(\beps_0) \right\|_{\Le_2} 
		&\leq \sum_{t=2}^n |g_{t,n}| \left\| G_{t,n} - G_{t-1,n} \right\|_{\Le_2} + \|G_{t-1,n}\|_{\Le_2} |g_{t,n} - g_{t-1,n}|  \\
		&\leq \Phi_n \Theta_n \Gamma_n + \Psi_n \Theta_n.
	\end{align*}
	We now apply Theorem \ref{thm:Gauss-ts}:
	Upon potentially enlarging the probability space, there exist independent Gaussian random vectors $Y_{t,n} \sim \mathcal{N}(0, \Sigma_{t,n})$ such that
	\begin{align*}
		\left\| \max_{k=1,\ldots, n} \frac{1}{\sqrt{n}} \sum_{t=1}^k Z_{t,n} - Y_{t,n} \right\|_{\Le_2} 
		&\leq C_{\beta,q} \Phi_n \Theta_n  \left( \Phi_n\Theta_n \Gamma_n + \Psi_n\Theta_n \right)^{\frac{\beta-1}{2\beta}} \sqrt{\log(n)} \left( \tfrac{m}{n} \right)^{\xi(q,\beta)}.
	\end{align*}

	\subsection*{Proof of Lemma \ref{lem:FCLT}}
	Define the processes
	\begin{align*}
		I_n(u) = \int_0^u \Sigma_{\lfloor vn\rfloor,n}^\frac{1}{2}\, dW_v,\qquad
		J_n(u) = \frac{1}{\sqrt{n}} \sum_{t=1}^{\lfloor un\rfloor} Y_{t,n},
	\end{align*}
	such that $(J_n(u))_{u\in[0,1]} \deq (I_n(\lfloor un\rfloor /n))_{u\in[0,1]}$. 
	Because $\Sigma_{t,n}$ is bounded, we find that  $\sup_{u\in[0,1]}|I_n(\lfloor un\rfloor /n) - I_n(u)| \pconv 0$.
	The process $I_n(u)$ is fully determined via its quadratic variation $[I_n]_u = \int_0^u \Sigma_{\lfloor vn\rfloor,n}\, dv$.
	To establish the weak convergence, it suffices to show that $[I_n]_u \to \int_0^u \Sigma_v\, dv$ for each $u$.
	Note that this convergence will then be uniform by monotonicity of the limit, and because the $\Sigma_{t,n}$ are bounded.
	
	To show convergence of $[I_n]_u$, we denote $\gamma_{t,n}(h) = \Cov[G_{t,n}(\beps_0), G_{t,n}(\beps_h)]$ and $\gamma_{u}(h) = \Cov[G_{u}(\beps_0), G_{u}(\beps_h)]$, and observe that
	\begin{align*}
		\|\Sigma_{t,n} - \Sigma_v\| 
		&\leq \sum_{h=-\infty}^\infty \|g_{t,n}-g_v\| \|\gamma_{t,n}(h)\| \|g_{t,n}\| \\
		&\quad + \sum_{h=-\infty}^\infty \|g_v\| \|\gamma_{t,n}(h)-\gamma_v(h)\| \|g_{t,n}\| \\
		&\quad + \sum_{h=-\infty}^\infty \|g_v\| \|\gamma_{t,n}(h)\| \|g_{t,n}-g_v\| \\
		&\leq C\|g_{t,n} - g_v\| \sum_{h=-\infty}^\infty (|h|+1)^{-\beta} +  C\sum_{h=-\infty}^\infty \|G_{t,n} - G_v\|_{\Le_2} \wedge (|h|+1)^{-\beta}.
	\end{align*}
	Thus,
	\begin{align*}
		&\int_0^1 \|\Sigma_{\lfloor vn\rfloor,n} - \Sigma_v\|\, dv \\
		&\leq C \int_0^1 \|g_{\lfloor vn\rfloor,n} - g_v\|\, dv + \sum_{h=-\infty}^\infty \left[ \int_0^1 \|G_{\lfloor vn\rfloor,n} - G_v\|_{\Le_2}\, dv \wedge (|h|+1)^{-\beta}\right].
	\end{align*}
	By virtue of assumptions \eqref{ass:A3} and \eqref{ass:A4}, and by dominated convergence, the latter term tends to zero as $n\to\infty$, thus completing the proof.

	\subsection*{Proof of Theorem \ref{thm:studentized}}
	Because of the Lipschitz-continuity of $\Sigma$, $X_{t,n}$ satisfies \eqref{ass:A1} and \eqref{ass:A2} for any $\beta>1$, with $\Gamma_n=1$ and $\Theta_n = \Theta(\beta) = \mathcal{O}(1)$. 
	Assumption \eqref{ass:A3} also holds, with $G_u(\beps_t) = \Sigma(u)^{\frac{1}{2}} \eta_t$.
	
	The multipliers are given by $\widehat{g}_{t,n} = \widehat{\Sigma}(t/n)^{-\frac{1}{2}}$, and we define $g_u = \Sigma(u)^{-\frac{1}{2}}$ and
	\begin{align*}
		\begin{cases}
			\Sigma(t/n)^{-\frac{1}{2}}, & t = k_n+1,\ldots, n, \\
			c I_{d\times d}, & t=1,\ldots, k_n,
		\end{cases} 
	\end{align*}
	Then $\Psi_n,\Phi_n = \mathcal{O}(1)$, \eqref{ass:A4} holds, and \eqref{ass:AX} with $L_n=1$ and $\Theta_n=0$, due to independence of the $X_{t,n}$. 
	
	It remains to compute $\Lambda_n$. 
	A standard bias-variance decomposition exploiting the Lipschitz continuity of $\Sigma(u)$ yields, for $t=k_n+1,\ldots, n$,
	\begin{align*}
		\EE \|\widetilde{\Sigma}(t/n) - \Sigma(t/n)\|^2 \;\leq\; C \left( \left(\frac{k_n}{n}\right)^2 + \frac{1}{k_n} \right) 
		\;\leq\; C n^{-\frac{2}{3}},
	\end{align*}
	for $k_n \asymp n^{-\frac{1}{3}}$ and a constant $C$ changing from line to line, but depending on neither $n$ nor $t$.
	Because the mapping $A\mapsto A^{-\frac{1}{2}}$ is Lipschitz-continuous on the set $\{A\succeq c I_{d\times d}\}$, we readily find that $\EE \|\widehat{g}_{t,n}-g_{t,n}\|^2 \leq C n^{-\frac{2}{3}}$ for all $t=1,\ldots, n$. 
	Thus, $\Lambda_n = \mathcal{O}(n^{\frac{1}{6}})$, which shows that \ref{thm:multiplier-FCLT} is applicable.

	\bibliography{strong-approximation.bib}
	\bibliographystyle{apalike}
	
\end{document}